\documentclass[12pt]{amsart}
\usepackage{amssymb}
\usepackage[latin1]{inputenc}
\usepackage{amsthm}
\usepackage{amsfonts}
\usepackage{mathrsfs}
\usepackage{graphicx}
\usepackage{amsmath}
\let\mathcal\mathscr
\def\bB{{\mathbb B}}
\def\bC{{\mathbb C}}

\def\bP{{\mathbb P}}

\newtheorem{theorem}{Theorem}

\newtheorem{conjecture}[theorem]{Conjecture}
\newtheorem{corollary}[theorem]{Corollary}

\newtheorem{definition}[theorem]{Definition}
\newtheorem{example}[theorem]{Example}

\newtheorem{proposition}[theorem]{Proposition}
\newtheorem{prop}[theorem]{Proposition}

\begin{document}
\thanks{{ E. R.
partially supported by a CIRGET fellowship and by the Chaire de
Recherche du Canada en alg\`ebre, combinatoire et informatique
math\'ematique de l'UQAM}}

\title{On the logarithmic Kobayashi conjecture}
\author{Gianluca Pacienza and Erwan Rousseau}
\date{}
\maketitle

%
\begin{quote}
{\small{\bf Abstract.} We study the hyperbolicity of the log
variety $(\mathbb{P}^n, X)$, where $X$ is a very general
hypersurface of degree $d\geq 2n+1$ (which is the bound predicted
by the Kobayashi conjecture). Using a positivity result for the
sheaf of (twisted) logarithmic vector fields, which may be of
independent interest, we show that any log-subvariety of
$(\mathbb{P}^n, X)$ is of log-general type, give a new proof of
the algebraic hyperbolicity of $(\mathbb{P}^n, X)$, and exclude
the existence of maximal rank families of entire curves in the
complement of the universal degree $d$ hypersurface. Moreover, we
prove that, as in the compact case, the algebraic hyperbolicity of
a log-variety  is a necessary condition for the metric one. }
\end{quote}

%
\section{Introduction}
%

A complex manifold $V$ is {\it hyperbolic} in the sense of S.
Kobayashi if the hyperbolic pseudodistance defined on $V$ (see
section 2 for precise definitions) is a distance.
 A necessary condition for the hyperbolicity of a
complex manifold is the constancy of holomorphic maps from
$\mathbb C$ to $V$. Hyperbolic complex manifolds have been studied
in the two following contexts. One is the hyperbolicity of a
compact complex manifold, in which case, thanks to a criterion due
to R. Brody \cite{Bro}, the above necessary condition is also
sufficient to guarantee the hyperbolicity of $V$. The other is the
hyperbolicity of a compact complex manifold with an ample divisor
removed. In the case of complements of projective hypersurfaces we
have the Kobayashi conjecture \cite{Ko70}:

%

\begin{conjecture}
\label{c2} The complement $\mathbb{P}^{n}\backslash X$   of a
general hypersurface $X\subset \mathbb{P}^{n}$ of degree $\deg
X\geq 2n+1,\ n\geq 2,$ is hyperbolic.
\end{conjecture}

Notice that the lower bound in Conjecture \ref{c2} is sharp,
since, as noticed first by M. Za\u\i denberg \cite{Z}, there
exists a line intersecting a general degree $2n$ hyperpersurface
in two points. (For the Kobayashi conjecture on the hyperbolicity
of a general hypersurface of high degree, and the related results,
we refer the reader to \S 2.2).


In the present paper we study questions related to Conjecture
\ref{c2} (which is proved for $n=2$ and $d\geq 15$ in \cite{E.G}),
by extending to the logarithmic setting (part of) the techniques
and ideas successfully used  in the compact case.

Let $\overline{V}$ be a variety with a normal crossing divisor
$D.$ The pair $(\overline{V},D)$ is called a log-variety. Let
$V=\overline{V}\backslash D.$ We denote by $\overline{T}_{V}^{\ast
}=T_{\overline{V}}^{\ast }(\log D)$ its log-cotangent bundle and
by $\overline{K}_{V}=\wedge ^{\dim(V)}\ \overline{T}_{V}^{\ast }=
K_{\overline{V}%
}(D)$ its log-canonical bundle. In the third section, in order to
study the algebraic hyperbolicity properties of the log-variety
$(\mathbb{P}^n, X)$, where $X$ is a very general hypersurface, we
prove the following non-vanishing result.

\begin{theorem}
\label{t1}Let $X$ be a very general hypersurface of arbitrary
degree $d$
 in $\mathbb{P}^{n}$. Let $Y$ be a $k$-dimensional subvariety
in $\mathbb{P}^{n}$ meeting X properly, $D:=Y\cap X$ the induced
divisor and $\nu: \widetilde{Y}\to Y$ a log-resolution of $(Y,D)$
i.e $\widetilde{Y}$ is smooth, $\nu$ is a projective birational
morphism and $\nu ^{-1}(D)+Exc(\nu)$ is a normal crossing divisor.
Then
$$
 h^0(\widetilde Y, \overline{K}_{\widetilde Y}\otimes \nu^*
\mathcal{O}_{\mathbb{P}^n}(2n+1-k-d))\not= 0,
$$
where $\overline{K}_{\widetilde Y}$ denotes the log-canonical
bundle of the log-variety $(\widetilde Y,\nu^{-1} (D))$.
\end{theorem}

In particular we deduce :

\begin{corollary}\label{cor1}
Let $X\subset \mathbb{P}^{n}$ be a very general hypersurface of
degree $d\geq 2n+2-k,\ k\geq1$. Then any $k$-dimensional
log-subvariety $(Y,D)$ of $(\mathbb{P}^n,X)$, for $Y$ not
contained in $X$, is of log-general type, that is, any
log-resolution $\nu:\widetilde Y\to Y$ of $(Y,D)$ has big
log-canonical bundle ${K}_{\widetilde Y} (\nu^{-1}(D))$.
\end{corollary}

\begin{corollary}\label{cor2}
Let $X\subset \mathbb{P}^{n}$ be a very general hypersurface of
degree $d\geq 2n+1$, and $C\subset \mathbb{P}^n$ a curve not
contained in $X$. Then
$$
2g(\widetilde{C})-2+i(C,X)\geq (d-2n)\deg C
$$
where $\nu: \widetilde{C}\to C$ is the normalization of $C$,
$g(\widetilde{C})$ its genus, and $i(C,X)$ is the number of
distinct points in $\nu^{-1}(X)$.

\end{corollary}

The inequality in Corollary \ref{cor2} has been previously proved
by Xi Chen in \cite{ch01} by means of a delicate degeneration
argument (see also \cite{V2}, Theorem 3.10, where it is proved
that any $C$ intersects a general hypersurface of degree
$d=2n-2+r,\ r\geq 3,$ in at least $r$ points, as well as
\cite{ch00} and \cite{Xu} for the case $n=2$). Note that, as a
consequence, one gets that there is no entire curve
$f:\mathbb{C}\to \mathbb{P}^n \setminus X$ in the complement of a
very general hypersurface of degree $d\geq 2n+1,\ n\geq2$, if the
Zariski closure $\overline{f(\mathbb{C})}$ is an algebraic curve.
Notice moreover that both Corollaries \ref{cor1} and \ref{cor2}
are sharp: the latter by the result of Za\u \i denberg we have
quoted above, and the former by a natural generalization of Za\u
\i denberg's result, which we present at the end of \S 3.

In the fourth section we prove that, as one expects, the
hyperbolicity of a log-variety implies its algebraic
hyperbolicity, thus answering a question raised by Xi Chen in
\cite{ch01}.

\begin{theorem}
\label{t2} \textit{Let }$\mathit{X}$\textit{\ be a projective manifold and }$%
\mathit{D}$\textit{\ an effective divisor on }$\mathit{X}$ such
that $\mathit{X\backslash D}$ is hyperbolic and hyperbolically
imbedded. Let $\omega $ be a hermitian metric on $X$. Then there
exists $\varepsilon
>0$ such that
\begin{equation*}
2g(\widetilde{C})-2+i(C,D)\geq \varepsilon \deg _{\omega }(C)
\end{equation*}
for every compact irreducible curve $C\subset X$ with $C\nsubseteq
D,$ where
$\widetilde{C}$ is the normalization of $C,$ $g(\widetilde{C})$ its genus and $%
\deg _{\omega }(C)=\int_{C}\omega .$
\end{theorem}

In the last section we strenghten the conclusion of Corollary
\ref{cor2}, and prove that there is no entire curve, varying in a
family of maximal rank, in the complement of a general
hypersurface of degree at least $2n+1$, without assumptions on the
Zariski closure of the entire curve.

\begin{theorem}\label{maxrank}
Let $U\to
\bP^{N_d}:=\mathbb{P}H^0(\bP^n,\mathcal{O}_{\mathbb{P}^n}(d))$ be
an \'etale cover of an open subset of $\bP^{N_d}$, and let
$\Phi:\bC\times U\to \bP^{n}\times U$ be a holomorphic map such
that $\Phi(\bC\times\{ t\})\subset \bP^{n} \setminus X_t$ for all
$t\in U$. If $d\ge 2n+1 $,   the rank of $\Phi$ cannot be maximal
anywhere.
\end{theorem}
The above theorem is of course a consequence of Conjecture
\ref{c2}, and represents the logarithmic analog of the main result
of \cite{DPP}. Again, by Za\u \i denberg's example, our result is
sharp, as for $d=2n$ one can consider the family of the
exponential maps associated to the family of lines $\ell_t$
intersecting the hypersurface $X_t$ in two points :
$$
 \exp_t : \mathbb{C}\times \{t\}\to \mathbb{C}^*=
 \ell_t\setminus (\ell_t \cap X_t)\subset \mathbb{P}^n\setminus X_t.
$$

We denote by $\mathcal{X}\subset \mathbb{P}^n\times
\mathbb{P}^{N_d}$ the family of degree $d$ hypersurfaces. The
proofs of Theorems \ref{t1} and \ref{maxrank} use the global
generation of the sheaf of (twisted) vector fields on the log
variety $(\mathbb{P}^n\times \mathbb{P}^{N_d}, \mathcal{X})$ (see
\S 3, Proposition \ref{prop1} for the precise statement). Once the
logarithmic framework is set, our approach allows natural proofs,
which are formally equal to those of the corresponding
hyperbolicity properties of $X\subset \mathbb{P}^n$. In that
sense, it unifies the compact and the logarithmic cases. Notice
moreover that the analogous global generation result in the
compact case is the first step in Y.-T. Siu's proof of the
hyperbolicity of a very general hypersurface $X\subset
\mathbb{P}^n$, for $d\gg n$ (see \cite{SY04}, Lemma 4, and
\cite{V'}, Proposition 1.1). It seems then plausible that using  a
generalization of Proposition \ref{prop1} to logarithmic jet
bundles and following the strategy outlined in \cite{SY04}, one
could prove Conjecture \ref{c2} for very high degree : this has
been done for $n=3$ in  \cite{Rou06bis}.



%
\section{Preliminaries}
%
\subsection{Log-manifolds}
Let $\overline{V}$ be a complex manifold with a normal crossing
divisor $D.$ The pair $(\overline{V},D)$ is called a {\it
log-manifold}. Let $V=\overline{V}\backslash D$ be the complement
of $D$.

Following \cite{Ii}, the logarithmic
cotangent sheaf $\overline{T}_{V}^{\ast }=T_{\overline{V}%
}^{\ast }(\log D)$ is defined as the locally free subsheaf of the
sheaf of
meromorphic 1-forms on $\overline{V},$ whose restriction to $V$ is $%
T_{V}^{\ast }$ and whose localization at any point $x\in D$ is
given by
\begin{equation*}
\overline{T}_{V,x}^{\ast }=\underset{i=1}{\overset{l}{\sum }}\mathcal{O}_{%
\overline{V},x}\frac{dz_{i}}{z_{i}}+\underset{j=1+1}{\overset{n}{\sum }}%
\mathcal{O}_{\overline{V},x}dz_{j}
\end{equation*}
where the local coordinates $z_{1,}...,z_{n}$ around $x$ are
chosen such that $D=\{$ $z_{1}...z_{l}=0\}.$

Its dual, the logarithmic tangent sheaf $\overline{T}_{V}=T_{\overline{V}%
}(-\log D)$ is a locally free subsheaf of the holomorphic tangent bundle $T_{%
\overline{V}},$ whose restriction to $V$ is $T_{V}$ and whose
localization at any point $x\in D$ is given by
\begin{equation*}
\overline{T}_{V,x}=\underset{i=1}{\overset{l}{\sum }}\mathcal{O}_{\overline{V%
},x}z_{i}\frac{\partial }{\partial z_{i}}+\underset{j=1+1}{\overset{n}{\sum }%
}\mathcal{O}_{\overline{V},x}\frac{\partial }{\partial z_{j}}.
\end{equation*}

Recall that starting with an arbitrary divisor, Hironaka's theorem
on resolution of singularities \cite{Hi} guarantees that we can
replace it by a normal crossing one after performing some
blowing-ups.

\begin{theorem}
Let $V$ be an irreducible complex algebraic variety (possibly
singular), and let $D\subset V$ be an effective Cartier divisor on
$V.$ There is a projective birational morphism
\begin{equation*}
\mu :V^{\prime }\rightarrow V,
\end{equation*}
where $V^{\prime }$ is non singular and $\mu $ has divisorial
exceptional locus $except(\mu ),$ such that
\begin{equation*}
\mu ^{-1}(D)+except(\mu )
\end{equation*}
is a normal crossing divisor.
\end{theorem}

One calls $V^{\prime}$ a {\it log-resolution} of $(V,D)$.

\subsection{Hyperbolicity and algebraic hyperbolicity}
Let $X$ be a complex manifold. We denote by $f:\Delta \rightarrow
X$ an
arbitrary holomorphic map from the unit disk $\Delta \subset \mathbb{C}$ to $%
X.$ The Kobayashi-Royden infinitesimal pseudo-metric \cite{Ko70}
on $X$ is the Finsler pseudometric on the tangent bundle $T_{X}$
defined by
\begin{equation*}
k_{X}(\xi )=\inf \{\lambda >0;\exists f:\Delta \rightarrow
X,f(0)=x,\lambda f^{\prime }(0)=\xi \},\text{ }x\in X,\xi \in
T_{X,x}.
\end{equation*}
The Kobayashi pseudo-distance $d_{X}$, is the geodesic
pseudodistance obtained by integrating the Kobayashi-Royden
infinitesimal pseudometric. The manifold $X$ is {\it hyperbolic}
in the sense of S. Kobayashi if the hyperbolic pseudodistance
defined on $X$ is a distance.

Directly from the definition of the Kobayashi pseudo-distance one
can see that if $f:X\rightarrow Y$ is a holomorphic map of complex
manifolds then it is distance decreasing i.e for $x,x^{\prime }\in
X$ we have
\begin{equation*}
d_{Y}(f(x),f(x^{\prime }))\leq d_{X}(x,x^{\prime }).
\end{equation*}

As mentioned in the introduction, in the case of a general
projective hypersurface $X$, both $X$ and its  complement are
conjectured to be hyperbolic, as soon as the degree of $X$ is high
enough.

\begin{conjecture}
\label{con1}A general hypersurface $X\subset \mathbb{P}^{n},\
n\geq 3,$ of degree $\deg X\geq 2n-1$ is hyperbolic.
\end{conjecture}

(For $n\geq 4$, the natural lower bound should be $2n-2$).

The most  important confirmation of Conjecture \ref{con1} has been
obtained by Y.-T. Siu \cite {SY04}, who proved it for  $d_{n}\gg
n$. As for the known lower bounds on the degree, Conjecture
\ref{con1} has been studied  for $n=3$ in \cite{DEG00} and
\cite{MQ}, where the bound $d\geq 21$ (respectively $d\geq 36$)
has been obtained (in the recent preprint \cite{Pau}, this bound
has been improved to $d\geq 18$). In \cite{Rou06}, the second
author proved a weak form of Conjecture \ref{con1} for $n=4$ and
$d\geq 593.$

It is widely believed that when dealing with a projective variety
$V$, there should exist a property of algebraic nature equivalent
to the hyperbolicity of $V$. In the compact case, Demailly (see
\cite{De95}) introduced the notion of {\it algebraic
hyperbolicity}, and proved it is a necessary condition for the
hyperbolicity. S. Lang proposed another property, namely the fact
that any subvariety of $V$ is of general type. Both properties
have been checked for very general hypersurfaces of degree $d\geq
2n-2$, in \cite{CR} and \cite{P}, building on ideas and techniques
introduced by C. Voisin \cite{V} (see also \cite{Cle} and
\cite{Ein}).

In analogy to the compact case, Xi Chen \cite{ch01} studied the
notion of {algebraic hyperbolicity}, in the sense of Demailly, for
log-manifolds.

\begin{definition}
{\em Let $(X,D)$ be a log-manifold. For each reduced curve
$C\subset X$ that
meets $D$ properly, let $\nu :%
\widetilde{C}\rightarrow C$ be the normalization of $C.$ Then
$i(C,D)$ is the number of distinct points in the set $\nu
^{-1}(D)\subset \widetilde{C}.$}
\end{definition}

\begin{definition}\label{alghyp}
{\em A logarithmic variety $(X,D)$ is algebraically hyperbolic if
there exists a positive number $\varepsilon $ such that
\begin{equation*}
2g(\widetilde{C})-2+i(C,D)\geq \varepsilon \deg _{\omega }(C)
\end{equation*}
for all reduced and irreducible curves $C\subset X$ meeting $D$
properly where
$\widetilde{C}$ is the normalization of $C,$ $g(\widetilde{C})$ its genus and $%
\deg _{\omega }(C)=\int_{C}\omega$ with $\omega$ a hermitian
metric on $X$.}
\end{definition}

In the next section, we prove that the algebraic hyperbolicity in
the sense of Demailly, as well as the algebraic property analogous
to that proposed by Lang, hold for the complement of a very
general projective hypersurface of degree at least equal to
$2n+1$.

%
\section{Algebraic hyperbolicity of the log variety $(\mathbb{P}^n,X)$}
%

In this section we give the proof of theorem \ref{t1}, using
logarithmic techniques, and the global generation of the sheaf of
(twisted) logarithmic vector fields, which we now introduce.


\smallskip

Fix the following notations:

$\mathbb{P}^{N_d}:=
\mathbb{P}H^{0}(\mathbb{P}^{n},\mathcal{O}_{\mathbb{P}^{n}}(d))$
denotes the parameter space for degree $d$ hypersurfaces in
$\mathbb{P}^n$.

$\mathcal{X}\subset \mathbb{P}^{n}\times \mathbb{P}^{N_d}$ denotes
the universal hypersurface of degree $d$, and $p$ and $q$ the
projections of $\mathbb{P}^{n}\times \mathbb{P}^{N_d}$ onto the
two factors.

$X_{F}\subset \mathbb{P}_{F}^{n}$ is the hypersurface defined by
the homogeneous polynomial $F\in \mathbb{P}^{N_d}.$

For a smooth hypersurface $X_{F}$ we have the corresponding
logarithmic
manifold $(\mathbb{P}_{F}^{n},X_{F})$, with logarithmic tangent sheaf $%
\overline{T}_{\mathbb{P}_{F}^{n}}=T_{\mathbb{P}_{F}^{n}}(-\log
X_{F}),$ logarithmic cotangent sheaf $\overline{\Omega
}_{\mathbb{P}_{F}^{n}}=\Omega
_{\mathbb{P}_{F}^{n}}(\log X_{F})$ and logarithmic canonical sheaf $%
\overline{K}_{\mathbb{P}_{F}^{n}}=K_{\mathbb{P}_{F}^{n}}\otimes \mathcal{O}%
(d)=\mathcal{O}(d-n-1).$

\bigskip

\subsection{The global generation result}
We shall extend to the logarithmic setting an approach initiated
by Clemens, Ein, Voisin and Siu (see \cite{Cle}, \cite{Ein}, \cite
{V}, \cite{SY04}):

\begin{proposition}
\label{prop1}The twisted logarithmic tangent bundle
$$
 T_{\mathbb{P}^{n}\times
 \mathbb{P}^{N_d}}(-\log \mathcal{X)}(1,0):=
 T_{\mathbb{P}^{n}\times
 \mathbb{P}^{N_d}}(-\log \mathcal{X)}\otimes p^*\mathcal{O}_{\bP^n}(1)
$$
is generated by its global sections.
\end{proposition}

\begin{proof}[Proof]
Let $\mathcal{X}\subset \mathbb{P}^{n}\times \mathbb{P}^{N_{d}}$
be the universal hypersurface of degree $d$ given by the equation
\begin{equation*}
\underset{\left| \alpha \right| =d}{\sum }a_{\alpha }Z^{\alpha }=0
\end{equation*}
where $[a]\in \mathbb{P}^{N_{d}}$ and $[Z]\in \mathbb{P}^{n},$ for
$\alpha =(\alpha _{0},...,\alpha _{n})\in \mathbb{N}^{n+1},$
$\left| \alpha \right| =\sum_{i}\alpha _{i}$ and if
$Z=(Z_{0},...,Z_{n})$ are homogeneous coordinates on
$\mathbb{P}^{n},$ then $Z^{\alpha }=\prod Z_{j}^{\alpha _{j}}.$
Notice that $\mathcal{X}$ is a smooth hypersurface of bidegree
$(d,1)$
in $\mathbb{P}^{n}\times \mathbb{P}^{N_{d}}.$ We consider the log-manifold $(%
\mathbb{P}^{n}\times \mathbb{P}^{N_{d}},\mathcal{X)}$. Let us
consider
\begin{equation*}
\mathcal{Z}=(a_{0...0d}Z_{n+1}^{d}+\underset{\left| \alpha \right|
=d,\alpha _{n+2}=0}{\sum }a_{\alpha }Z^{\alpha }=0)\subset
\mathbb{P}^{n+1}\times U
\end{equation*}
where $\alpha \in \mathbb{N}^{n+2},$ and
$$
 U:=(a_{0...0d}\neq 0)\cap
 \left( \underset{\left| \alpha \right| =d,\alpha _{n+2}=0}{\cup
 }(a_{\alpha }\neq
 0)\right) \subset \mathbb{P}^{N_{d}+1}.
$$
Consider the natural projection $\pi :\mathcal{Z}\rightarrow
\mathbb{P}^{n}\times \mathbb{P}^{N_{d}}$ and set
$$
 \mathcal{H}:=\pi^{-1}(\mathcal{X})=\{Z_{n+1}=0\}.
$$
Therefore we obtain a dominant log-morphism $%
\pi :(\mathcal{Z},\mathcal{H})\rightarrow (\mathbb{P}^{n}\times
\mathbb{P}^{N_{d}},\mathcal{X})$ which induces a
 map
\begin{equation*}
\pi _{\ast
}:\overline{T}_{\mathcal{Z}}(1,0):=T_{\mathcal{Z}}(-\log
\mathcal{H})(1,0)\rightarrow \overline{T}_{\mathbb{P}^{n}\times \mathbb{P}%
^{N_{d}}}(1,0):=T_{\mathbb{P}^{n}\times \mathbb{P}^{N}}(-\log
\mathcal{X)}(1,0).
\end{equation*}

Therefore we are reduced to prove that
$\overline{T}_{\mathcal{Z}}(1,0)$ is generated by its global
sections. Consider the open set $U_{0}=\{Z_{0}\neq 0\}\times U$ in
$\mathbb{P}^{n+1}\times U$ with the induced inhomogeneous
coordinates. The equation of $\mathcal{Z}$ on $U_{0}$ becomes
\begin{equation*}
\mathcal{Z}_{0}:=\{z_{n+1}^{d}+\underset{\left| \alpha \right|
\leq d,\alpha _{n+1}=0}{\sum }a_{\alpha }z^{\alpha }=0\}.
\end{equation*}

Consider the vector field
\begin{equation*}
V_{\alpha ,j}=\frac{\partial }{\partial a_{\alpha }}-z_{j}\frac{\partial }{%
\partial a_{\widetilde{\alpha }}}
\end{equation*}
where $\alpha \in \mathbb{N}^{n+1},\alpha _{n+1}=0,\exists j$
$\alpha
_{j}\geq 1,\widetilde{\alpha _{k}}=\alpha _{k}$ if $k\neq j$ and $\widetilde{%
\alpha _{j}}=\alpha _{j}-1.$ Notice that $V_{\alpha ,j}$ is a
logarithmic vector field of
$(\mathcal{Z}_{0},\mathcal{H}_{0}:=\mathcal{H}\cap U_0)$ which
extends to $(\mathcal{Z},\mathcal{H})$ with a pole order equal to
1.

Consider a vector field
\begin{equation*}
V_{0}=\overset{n}{\underset{j=1}{\sum }}v_{j}\frac{\partial }{\partial z_{j}}%
+v_{n+1}z_{n+1}\frac{\partial }{\partial z_{n+1}}
\end{equation*}
where
$$
 v_{j}=\overset{n}{\underset{k=1}{\sum }}v_{k}^{(j)}z_{k}+v_{0}^{j},1%
 \leq j\leq n
$$
is linear in the variables $z_k$, and $v_{n+1}\in \mathbb{C}.$ We
claim that there exists a logarithmic vector field
\begin{equation*}
V=\underset{\left| \alpha \right| \leq d,\alpha _{n+1}=0}{\sum }v_{\alpha }%
\frac{\partial }{\partial a_{\alpha }}+V_{0}
\end{equation*}
tangent to $\mathcal{Z}_{0}.$ Indeed, the condition to be
satisfied is
\begin{equation*}
\underset{\alpha }{\sum }v_{\alpha }z^{\alpha }+\underset{\alpha ,j}{\sum }%
a_{\alpha }v_{j}\frac{\partial z^{\alpha }}{\partial z_{j}}=0
\end{equation*}
and the complex numbers $v_{\alpha }$ are chosen such that the
coefficient of $z^{\alpha }$ in the above equation is equal to
zero. This logarithmic vector field of
$(\mathcal{Z}_{0},\mathcal{H}_{0})$ extends to
$(\mathcal{Z},\mathcal{H}).$

The previous vector fields give the global generation of $\overline{T}_{%
\mathcal{Z}}(1,0).$
\end{proof}

\bigskip

\subsection{Sharp algebraic hyperbolicity properties for $(\mathbb{P}^n,X)$}
Having recorded in the previous subsection the needed positivity
result, we can now prove Theorem \ref{t1}, together with its
corollaries. We also show, in Example \ref{example}, that our
results are sharp in the degree.

\begin{proof}[Proof of Theorem \ref{t1}]
Let $U\subset \mathbb{P}^{N}$ be the open subset parametrizing
smooth hypersurface. We want to study families of $k$-dimensional
irreducible subvarieties inside $\mathbb{P}^{n}\times U$,
intersecting properly the family of hypersurfaces. So, eventually
passing to an \'etale cover of $U$, we consider an irreducible
subvariety $\mathcal{Y}\subset \mathbb{P}^{n}\times
\mathbb{P}^{N}$ such that the projection map
$\mathcal{Y}\rightarrow U$  is dominant of relative dimension $k$,
and such that $\mathcal{Y}$ intersects properly $\mathcal{X}$ (and
so does its generic fiber $Y_F$ with $X_F$). Let
$\mathcal{D}\subset \mathcal{Y}$ the family of divisors induced by
the intersections $D_F:=Y_F\cap X_F$. Let
$\widetilde{\mathcal{Y}}\rightarrow \mathcal{Y}$ be a log
resolution of $(\mathcal{Y},\mathcal{D})$ i.e
$\widetilde{\mathcal{Y}}$ is smooth and
$\widetilde{\mathcal{D}}=\nu ^{-1}(\mathcal{D})$ is a normal
crossing divisor, and so is its general fiber
$\widetilde{D}_F\subset \widetilde{Y}_F$.

For general $F\in U$, and {\it arbitrary} degree $d$ we want to
produce a non zero element in $ H^{0}(\widetilde{Y}_{F},
\overline{K}_{\widetilde{Y}_{F}}(2n+1-k-d))$, where
$\overline{K}_{\widetilde{Y}_{F}}= K_{\widetilde{Y}_{F}}(
\widetilde{D}_{F}).$



We have:

\begin{equation}\label{adjunction}
\text{ }\overline{K}_{\widetilde{Y}_{F}}\simeq \overline{\Omega }_{%
\widetilde{\mathcal{Y}}\left| \widetilde{Y}_{F}\right. }^{N+k}
\end{equation}

Indeed:

$$
 \Omega _{\widetilde{\mathcal{Y}}}^{N+k}\otimes
 \mathcal{O}_{\widetilde{\mathcal{Y}}}
 (\widetilde{\mathcal D})_{\left| \widetilde{Y}_{F}\right. }=
 \Omega _{\widetilde{%
 \mathcal{Y}}\left| \widetilde{Y}_{F}\right. }^{N+k}\otimes
 \mathcal{O}_{\widetilde{Y}_{F}}(\widetilde{D}_{F}),
$$

and by the adjunction formula
$$
 K_{\widetilde{Y}_{F}}=\Omega _{\widetilde{%
 \mathcal{Y}}\left| \widetilde{Y}_{F}\right. }^{N+k},
$$
since the normal bundle of a fiber in a family is trivial.

Using standard linear algebra, we have:

\begin{equation}\label{duality}
\text{ }(\bigwedge^{n-k}T_{\mathbb{P}^{n}\times
\mathbb{P}^{N}}(-\log
\mathcal{X)}_{\left| \mathbb{P}_{F}^{n}\right. })\otimes \overline{K}_{%
\mathbb{P}_{F}^{n}}\simeq \Omega _{\mathbb{P}^{n}\times \mathbb{P}%
^{N}}^{N+k}(\log \mathcal{X)}_{\left| \mathbb{P}_{F}^{n}\right. }
\end{equation}

The generically surjective map $K_{\mathbb{P}^{n}\times \mathbb{P}%
^{N}}(\mathcal{X)}\rightarrow K_{\widetilde{\mathcal{Y}}}
(\widetilde{\mathcal{D}}%
)$ induces a map
\begin{equation*}
\Omega _{\mathbb{P}^{n}\times \mathbb{P}^{N}}^{N+k}(\log \mathcal{X)}%
_{\left| \mathbb{P}_{F}^{n}\right. }(2n+1-k-d)
\rightarrow \overline{\Omega }_{%
\widetilde{\mathcal{Y}}\left| \widetilde{Y}_{F}\right.
}^{N+k}(2n+1+k-d)
\end{equation*}

that is non zero for $F$ general in $U.$

Recalling that
$$
 \overline{K}_{\mathbb{P}_{F}^{n}}=
 \mathcal{O}_{\mathbb{P}_{F}^{n}}(d-n-1)
$$
and using (\ref{adjunction}) and (\ref{duality}), it is enough to
show that
\begin{equation*}
(\bigwedge^{n-k}T_{\mathbb{P}^{n}\times \mathbb{P}^{N}}(-\log \mathcal{X)}%
_{\left| \mathbb{P}_{F}^{n}\right. })\otimes \mathcal{O}_{\mathbb{P}%
_{F}^{n}}(n-k)
\end{equation*}
is globally generated. To conclude, we notice that
\begin{equation*}
(\bigwedge^{n-k}T_{\mathbb{P}^{n}\times \mathbb{P}^{N}}(-\log \mathcal{X)}%
_{\left| \mathbb{P}_{F}^{n}\right. })\otimes
\mathcal{O}_{\mathbb{P}_{F}^{n}}(n-k)=
\bigwedge^{n-k}(T_{\mathbb{P}^{n}\times \mathbb{P}%
^{N}}(-\log \mathcal{X)}_{\left| \mathbb{P}_{F}^{n}\right. }(1))
\end{equation*}
and invoke the global
generation of $T_{\mathbb{P}^{n}\times \mathbb{P}%
^{N}}(-\log \mathcal{X)}_{\left| \mathbb{P}_{F}^{n}\right. }(1)$,
that follows from Proposition \ref{prop1}.

Letting the family $\mathcal{Y}$ vary, that is, varying the
Hilbert polynomial, we obtain that for $F$ outside a countable
union of proper closed subvarieties of $U$, all the
$k$-dimensional subvarieties $Y$ intersecting properly $X_F$
verify
\begin{equation}\label{not0}
 h^0(\widetilde Y, \overline{K}_{\widetilde Y}\otimes \nu^*
 \mathcal{O}_{\mathbb{P}^n}(2n+1-k-d))\not= 0.
\end{equation}
\end{proof}

Let us now show how to deduce Corollaries \ref{cor1} and
\ref{cor2}.

\begin{proof}[Proof of Corollary \ref{cor1}]
By the above theorem, the logarithmic canonical bundle of
$({\widetilde Y},{\widetilde D})$ may be written as the sum of the
{\it effective} line bundle $\overline{K}_{\widetilde Y}\otimes
\nu^*
 \mathcal{O}_{\mathbb{P}^n}(2n+1-k-d)$ and the line bundle
$\nu^* \mathcal{O}_{\mathbb{P}^n}(d-(2n+1-k))$. The latter is {\it
big}, as soon as $d\geq 2n+2-k$, so the corollary is proved.
\end{proof}

\begin{proof}[Proof of Corollary \ref{cor2}]
If $C\subset \mathbb{P}^{n}$ is a curve intersecting properly the
general hypersurface $X_F$, $f: \widetilde{C}\rightarrow C$ its
desingularization, $D:=C\cap X_F$ the divisor given by the
intersection with the hypersurface, and $\widetilde{D}=f^{-1}(D)$,
then by (\ref{not0}) we have
\begin{equation*}
0\leq \deg (K_{\widetilde{C}}(\widetilde{D})\otimes f^{\ast }\mathcal{O}_{\mathbb{P}%
^{n}}(2n-d))=2g(C)-2+i(C,X)-(d-2n)\deg C,
\end{equation*}
and we are done.
\end{proof}

A further consequence is the following.

\begin{corollary}
For a very general hypersurface X of degree $d\geq 2n+1$ in
$\mathbb{P}^{n},$
$\mathbb{P}^{n}\backslash X$ does not contain any algebraic torus $\mathbb{C}%
^{\ast }.$ Therefore a holomorphic map $f:\mathbb{C}\rightarrow \mathbb{P}%
^{n}\backslash X$ is constant if $f(\mathbb{C)}$ is contained in
an algebraic curve.
\end{corollary}

We end the present section by discussing an example which
generalizes \cite{Z}, and shows that also Corollary \ref{cor1} is
sharp.

\begin{example}\label{example}
\em{Given a degree $d$ hypersurface $X\subset \bP^n$ and an
integer $r\geq 1$, consider the bicontact locus
$\Delta_{r,d-r,X}\subset X$ consisting of points $x\in X$ through
which passes a line $\ell$ such that  $\ell$ has contact at least
$r$ at $x$ and, if it is not contained in $X$, $\ell$ intersects
$X$ is at most another point $x'$. In other words, generically we
have
$$
 \ell \cap X = r\cdot x + (d-r) \cdot x'.
$$
If $X$ is general of degree $d\leq 2n$, then $\Delta_{r,d-r,X}$ is
non empty and of the expected dimension $2n-d$ (see \cite{V}, or
\cite{P}, \S 4, for the proof of this fact, and for a description
of (a desingularization of) $\Delta_{r,d-r,X}$ as the zero locus
of a section of a vector bundle). Hence, taking $d=2n$, we recover
the existence of a line intersecting the general degree $d$
hypersurface in two points, as first observed in \cite{Z}. Now
take $d=2n+1-k$. In this case the dimension of $\Delta_{r,d-r,X}$
equals $k-1$. Let $Y$ be the ($k$-dimensional) subvariety of
$\bP^n$ spanned by the lines parametrized by $\Delta_{r,d-r,X}$.
If its desingularization $\widetilde{Y}$ were of log-general type,
then by restriction to the general line $\ell$ in the ruling of
$Y$, we would get non constant sections of $H^0(\ell,
\overline{K}_{\ell}^{\otimes m})$. So we would get to a
contradiction since $\overline{K}_{\ell} = \mathcal{O}_{\ell}$, as
$\ell$ intersects $X$ at two points. }
\end{example}
%
\section{Hyperbolicity and algebraic hyperbolicty of log varieties}
%
In this section, we would like to answer a question raised by X.
Chen in \cite{ch01}:

\textit{Let }$\mathit{X}$\textit{\ be a projective manifold and }$\mathit{D}$%
\textit{\ an effective divisor on }$\mathit{X}$\textit{. Is it true that if }%
$\mathit{X\backslash D}$\textit{\ is hyperbolic and hyperbolically
imbedded then }$\mathit{(X,D)}$\textit{\ is algebraically
hyperbolic?}

The answer is positive. Namely we have theorem \ref{t2}.

To give a proof we need the following results.

First, we need a Gauss-Bonnet formula in the non compact case. We
follow the approach of \cite{Far} which we recall for the
convenience of the reader. Let $M$ be a Riemann surface. $M$ is
said to be of \textit{finite type} if there exists a compact
Riemann surface $M^{\prime }$ such that $M^{\prime }\backslash M$
consists of finitely many points. The genus of $M$ is defined as
the genus of $M^{\prime }$.

A puncture of $M$ is defined to be a domain $%
D_{0}\subset M$ conformally equivalent to $\{z\in
\mathbb{C};0<\left| z\right| <1\}.$ We will identify $z=0$ with
the puncture $D_{0}$.

Recall that a Kleinian group $G$ is a subgroup of $PGL_2$ whose
action on $\mathbb{P}^1 $ is discontinuous at some point and that
a Kleinian group is called Fuchsian if there is a disc invariant
under the action. Let $G$ be a Fuchsian group acting on the unit
disk $\Delta . $ Let $\{x_{1},x_{2},...,x_{n}\}$ be the set of
points of $\Delta /G$ that are either punctures or ramified points
of the projection $\pi :\Delta \rightarrow $ $\Delta /G$. Let $\nu
_{j}$ be the ramification index of $\pi ^{-1}(x_{j})$ and set $\nu
_{j}=\infty $ for punctures. Let us assume that $\Delta/G$ is of
finite type. If $\pi $ is ramified over finitely many points, then
we will say that $G$ is of \textit{finite type over }$\Delta.$ We
let $g$ be the genus of $\Delta/G.$ We can define the
characteristic of $G:$%
\begin{equation*}
\chi =2g-2+\underset{j=1}{\overset{n}{\sum }}\left( 1-\frac{1}{\nu _{j}}%
\right) .
\end{equation*}

We can project the Poincar\'{e} metric $\frac{4\left| dz\right| ^{2}}{%
(1-\left| z\right| ^{2})^{2}}$ on $\Delta /G$ which gives the
hyperbolic metric of constant curvature $-1.$ We have the
following theorems:

\begin{theorem}
(\cite{Far}, p.233) The area of $\Delta /G$ with respect to the
hyperbolic metric is finite and
\begin{equation*}
Area(\Delta /G)=2\pi \chi .
\end{equation*}
\end{theorem}

\begin{theorem}
(\cite{Far}, p.234) Let $M$ be a Riemann surface and
$\{x_{1},x_{2},...\}$ a discrete sequence on $M.$ To each point
$x_{k}$ we assign an integer $\nu _{k}\geq 2$ or $\infty .$

If $M=\mathbb{P}^1 $ we exclude two cases:

\begin{enumerate}
\item[(i)] $\{x_{1},x_{2},...\}$ consists of one point and $\nu
_{1}\neq \infty .$ \item[(ii)] $\{x_{1},x_{2},...\}$ consists of
two points and $\nu _{1}\neq \nu _{2}.$
\end{enumerate}

Let $M^{\prime }=M\backslash \underset{\nu _{k}=\infty }{\cup }%
\{x_{k}\}.$ Then there exists a simply connected Riemann surface
$\widetilde{M},$ a Kleinian group $G$ of self mappings of
$\widetilde{M}$ such that
\begin{enumerate}
\item[(a)] $\widetilde{M}/G\cong M^{\prime }$

\item[(b)] the natural projection $\pi :\widetilde{M}\rightarrow
M^{\prime }$ is
unramified except over the points $x_{k}$ with $\nu _{k}<\infty $ where the branch numbers verify $%
b_{\pi }(\widetilde{x})=\nu _{k}-1$ for all $\widetilde{x}\in \pi
^{-1}(\{x_{k}\}).$
\end{enumerate}
\end{theorem}

The third result we need is related to the notion of hyperbolic
imbedding (see \cite{Ko98}). Let $Z$ be a complex manifold and $Y$
a complex submanifold with compact closure $\overline{Y}.$ $Y$ is
hyperbolically imbedded in $Z$ if for every pair of distinct
points $p,q$ in $\overline{Y}\subset Z,$ there
exist neighborhoods $U_{p}$ and $U_{q}$ of $p$ and $q$ in $Z$ such that $%
d_{Y}(U_{p}\cap Y,U_{q}\cap Y)>0$ which is equivalent to say that $%
d_{Y}(p_{n},q_{n})$ cannot converge to zero when two sequences
$\{p_{n}\}$ and $\{q_{n}\}$ in $Y$ approach two distinct points
$p$ and $q$ of the boundary $\partial Y=\overline{Y}\backslash Y.$

Let us prove the following proposition which is another version of
Theorem 3.3.3 of \cite{Ko98}:

\begin{proposition}
Let $Y$ be a relatively compact complex submanifold (i.e
$\overline{Y}$ is compact) of a complex manifold $Z$. Then the
following are equivalent:
\begin{enumerate}
\item[(a)] $Y$ is hyperbolically imbedded in $Z$.

\item[(b)] Given a length function $L$ on $Z$ there is a positive constant $%
\varepsilon $ such that $k_{Y}\geq \varepsilon L$ on $Y.$
\end{enumerate}
\end{proposition}

\begin{proof}[Proof]
Let us prove that (a) implies (b). If $\varepsilon $ does not
exist then, from the definition of the Kobayashi infinitesimal
pseudometric, there exists a sequence $\{f_{n}\}$ of holomorphic
functions from $\Delta $ to $Y$ such that
\begin{equation*}
L(f_{n}^{\prime }(0))>n.
\end{equation*}

Since $\overline{Y}$ is compact we may assume that $\{f_{n}(0)\}$
converges to a point $p\in \overline{Y}.$

Let $U$ be a complete hyperbolic neighborhood of $p$ in $Z.$
Assume that there exists a positive number $r<1$ such that
$f_{n}(\Delta _{r})\subset U$ for $n\geq n_{0}.$ Then
$\{f_{n\left| \Delta _{r}\right. }:\Delta _{r}\rightarrow U\}$
would be relatively compact and woud have a subsequence which
converges to a holomorphic function from $\Delta _{r}$ to $U,$
which contradicts $L(f_{n}^{\prime }(0))>n.$

This means that for each positive integer $k,$ there exist a point
$z_{k}\in \Delta $ and an integer $n_{k}$ such that $\left|
z_{k}\right| <\frac{1}{k}$
and $f_{n_{k}}(z_{k})\notin U.$ Let $p_{k}=f_{n_{k}}(0)$ and $q_{k}=$ $%
f_{n_{k}}(z_{k}).$ By taking a subsequence we may assume that
$\{q_{k}\}$ converges to a point $q$ not in $U.$ Therefore we have
\begin{equation*}
d_{Y}(p_{k},q_{k})\leq d_{\Delta }(0,z_{k})\rightarrow 0\text{ for }%
k\rightarrow \infty ,
\end{equation*}

and this contradicts the fact that $Y$ is hyperbolically imbedded
in $Z$.

Let us prove that (b) implies (a). Let $\delta $ be the distance
function on $Z$ induced by $L$. Then
\begin{equation*}
\varepsilon \delta \leq d_{Y}\text{ on }Y
\end{equation*}

which implies obviously that $Y$ is hyperbolically imbedded in
$Z$.
\end{proof}

Now, we can prove Theorem \ref{t2}.

\begin{proof}[Proof of Theorem \ref{t2}]
Let $\nu :\widetilde{C}\rightarrow C$ be the normalization and $\widetilde{D}%
=\nu ^{-1}(D).$ As $X\backslash D$ is hyperbolic
$C^{\prime }=\widetilde{C}%
\backslash \widetilde{D}$ is hyperbolic and admits the unit disk
as its universal cover $\rho :\Delta \rightarrow C^{\prime }.$ Let
$k_{C^{\prime }}$ be the hyperbolic metric of constant curvature
$-1$ with $\mu _{C^{\prime }}$ its area element$.$ From the
distance decreasing property of Kobayashi metrics and the previous
proposition we have
\begin{equation*}
k_{C^{\prime }}(t)\geq k_{X\backslash D}(\nu _{\ast }(t))\geq
\varepsilon \left\| \nu _{\ast }(t)\right\| _{\omega },\text{
}\forall t\in T_{C^{\prime }}.
\end{equation*}

Therefore from the preceding two theorems we have
\begin{equation*}
2\pi (2g(\widetilde{C})-2+i(C,D))=\int_{C^{\prime }}\mu
_{C^{\prime }}\geq \varepsilon ^{2}\int_{C}\omega .
\end{equation*}
\end{proof}

%
\section{Families of entire curves in the complement of the universal hypersurface.}
%

The goal of this section is to prove that a family of entire
curves in the complement of the universal degree $d$ hypersurface
${\mathcal X}\subset \bP^n\times \bP^{N_d}$ cannot have maximal
rank, as soon as $d\geq 2n+1$ (as predicted by the log Kobayashi
conjecture). As in the compact case, which has been treated in
\cite{DPP},  such a result points out to the lack of a good
parameter space for entire curves. The proof goes exactly as in
the  compact case and relies on the global generation of the
bundle  $T_{\mathbb{P}^{n}\times \mathbb{P}^{N}}(-\log
\mathcal{X)}(1,0)$, proved in Proposition \ref{prop1}.

Let $U\to\bP^{N_d}$ be an \'etale cover of the open subset
parametrizing smooth hypersurfaces. As before, to render the
notation less heavy we will simply write $U\subset \bP^{N_d}$

Consider  a holomorphic map $\Phi: \bC\times U\to (\bP^n\times U)
\setminus {\mathcal X}$ over the base $U\subset \bP^{N_d}$.

As $U$ is an open set, we can shrink it and suppose that it is
equal to a polydisc $\bB(\delta_0)^{N_d}$. Consider the following
sequence of maps
$$
 \Phi_k: \bB(\delta_0k)^{N_d+ 1}\to (\bP^n\times U) \setminus{\mathcal X}
$$
 given by
$\displaystyle \Phi_k(z, \xi_1,\dots,\xi_{N_d})= \Phi (z k^{N_d},
\frac{1}{k}\xi_1,\dots,\frac{1}{k}\xi_{N_d})$. We will point out
where the change of the radius of the disc is used in the proof.

If $\Phi_1 =\Phi$ is of maximal rank, its Jacobian gives a nonzero
section
$$J_\Phi(z, \xi)=
{\frac{\partial \Phi }{\partial z}}\wedge {\frac{\partial \Phi
}{\partial \xi_1}}\wedge\dots \wedge {\frac{\partial \Phi
}{\partial \xi_{N_d}}}(z, \xi)\in \bigwedge ^{1+ N_d}
\Phi^*T_{{\mathcal X}, \Phi(z, \xi)}$$

Let us assume that $J_\Phi({\underline 0}) $ is nonzero in the
corresponding vector space. Remark   that, by construction,
$J_{\Phi_k}({\underline 0})= J_\Phi({\underline 0})$, for any
$k\geq 1$, hence $J_{\Phi_k}\in \Phi_k^*\Lambda ^{1+
N_d}T_{\bP^n\times \bP^{N_d}}$ is not identically zero.
Thanks to  Proposition \ref{prop1}, we can choose $(n-1)$
logarithmic vector fields
$$
 V_1,\dots,V_{n-1}\in H^0(T_{\mathbb{P}^{n}\times
\mathbb{P}^{N}}(-\log \mathcal{X)}(1,0) )$$ such that the sections
$$
 \sigma_k:= J_{\Phi_k}\wedge
 \Phi_k^*\bigl(V_1\wedge\dots\wedge V_{n-1}\bigr)\in
 \Phi_k^*(\overline{K}_{\bP^n\times\bP^{N_d}}^{\ -1}\otimes
 p^*{\mathcal O}_{\bP^n}(n-1))\bigr)
$$
are nonzero at the origin.

If $q$ is the projection of $\bP^n\times \bP^{N_d}$ on the
parameter space $\bP^{N_d}$, then under the assumption $d\geq
2n+1$, the restriction of
$\overline{K}_{\bP^n\times\bP^{N_d}}\otimes
 p^*{\mathcal O}_{\bP^n}(1-n)$ to
$q^{-1}(U)$ is ample (eventually after shrinking once again the
open subset $U$), hence we can endow this bundle with a metric $h$
of positive curvature.

For any $w\in \bB(\delta_0k)^{N_d+ 1}$ set
\begin{equation}\label{f_k}
 f_k (w)=
 \Vert \sigma_k (w) \Vert ^{2/(N_d+1)}_{\Phi_k ^*h^{-1}}.
\end{equation}

Notice that, by construction,  there exists a positive number $c$
such that for each $k\geq 1$, we have
\begin{equation}\label{f_k(0)}
   f_k(\underline{0})= c>0.
\end{equation}

On the other hand we have
\begin{prop}\label{tozero}
 For each $k\geq 1$ we have $f_k(\underline{0})\leq C\cdot k^{-2}$.
 In particular, as $k\to \infty$, we have $f_k(\underline{0})\to 0$.
\end{prop}
Theorem \ref{maxrank} follows from the fact that (\ref{f_k(0)})
and Proposition \ref{tozero} contradict each other.

We now give the proof of Proposition \ref{tozero}, which is very
close to that of the classical Ahlfors-Schwarz lemma.

\begin{proof}[Proof of Proposition \ref{tozero}]
First, notice that
 for each $k\geq 1$, there exists a positive constant $C$ such that
 we have
\begin{equation}\label{lowerbound}
 \Delta \log f_k\geq C\cdot f_k
\end{equation}

\noindent pointwise over the polydisc $\bB(\delta_0k)^{N_d+ 1}$.
Indeed, by construction, the image of the map $\Phi_k$ lies inside
$q^{-1}(U)$, for each $k\geq 1$, so that
$$
 i\partial{\bar \partial}
 \log\Vert \sigma_k \Vert ^2_{\Phi_k^*h^{-1}}\geq
 \Phi_k^*\Theta_h\bigl(\overline{K}_{\bP^n\times \bP^{N_d}}\otimes
 p^*{\mathcal O}_{\bP^n}(2-n)\bigr).
$$
Hence we get
\begin{eqnarray*}
 \Delta \log\Vert \sigma_k \Vert ^2_{\Phi_k^*h^{-1}}
 & \geq &
 C''\cdot
 \Bigl(\Big \Vert {\frac{\partial \Phi_k}{\partial z}}\Big \Vert ^2_\omega
 + \sum_{j=1}^{N_d}
\Big \Vert {\frac{\partial \Phi_k} {\partial \xi_j}}\Big \Vert ^2_\omega\Bigr)\\
& \geq &
 C'\cdot \Vert J_{\Phi_k}\Vert ^{2/(1+N_d)}_{\Lambda^{1+ N_d}\omega}\\
& \geq &
 C\cdot \Vert \sigma_k \Vert _{\Phi_k^*h^{-1}}^{2/(1+N_d)}
\end{eqnarray*}
and (\ref{lowerbound}) is proved (the above relations are obtained
using the vector inequalities
$$
 \Vert W_1 \wedge \dots \wedge W_s \Vert
 \leq \Vert W_1 \Vert \dots \Vert W_s \Vert
 \leq s^{-s} (\Vert W_1 \Vert +
 \dots + \Vert W_s \Vert)^s ).
$$

Then, consider the volume form of the Poincar\'e metric on the
polydisc
$$\psi_k= \frac{1}{\Bigl(1- \frac{\vert z\vert ^2}{\delta_0^2k^2}\Bigr)^2}
\prod_{j=1}^{N_d}\frac{1}{\Bigl(1- \frac{\vert \xi_j\vert
^2}{\delta_0^2k^2}\Bigr)^2}$$ A  computation shows that
\begin{equation}\label{upperbound}
\Delta \log \psi_k\leq C\cdot k^{-2}\psi_k.
\end{equation}
(Remark that the previous inequality can be obtained precisely
because, thanks to the reparameterization, we have the same radius
$\delta_0k$ for the components of the polydisc which is the domain
of $\psi_k$.)

Consider the function $\displaystyle (z, \xi)\mapsto \frac{f_k(z,
\xi)}{\psi_k(z, \xi)}$. Its maximum cannot be achieved at a
boundary point of the domain, since $\psi_k$ goes to infinity as
$(z, \xi)$ goes to the boundary. So at the maximum point $(z_0,
\xi_0)$, we have
$$\Delta \log f_k/\psi_k\leq 0.$$
 This inequality, combined with  (\ref{lowerbound})
and (\ref{upperbound}), gives
$$f_k(z_0, \xi_0)\leq C\cdot k^{-2}\psi_k(z_0, \xi_0)$$
Since the relation (5) is verified at the maximum point of the
quotient, the same is true at an arbitrary point, thus, in
particular, at the origin:
$$
f_k(\underline{0})\leq C \cdot k^{-2}.
$$

\end{proof}

\noindent
{\tt pacienza@math.u-strasbg.fr}\\
Institut de Recherche Math\'ematique Avanc\'ee\\
Universit\'e L. Pasteur et CNRS\\
7, rue Ren\'e Descartes, 67084 Strasbourg C\'edex - FRANCE

\medskip

\noindent
{\tt eroussea@math.uqam.ca}\\
D\'epartement de Math\'ematiques\\
Universit\'e du Qu\'ebec \`a Montr\'eal\\
Case Postale 8888, Succursale Centre-Ville\\
Montr\'eal (Qu\'ebec) CANADA H3C 3P8
\end{document}